\documentclass[a4paper,11pt,twoside]{article}
\date{Julho de 2014}

\usepackage{amssymb}
\usepackage[T1]{fontenc}
\usepackage[utf8]{inputenc}
\usepackage[english]{babel}
\usepackage{mathrsfs}
\usepackage{amsfonts,amsmath,lmodern,fancyhdr,lastpage,graphicx,nccfoots,caption}
\usepackage{afterpage}
\usepackage{amsthm}
\usepackage{enumerate}
\usepackage[all]{xy}
\usepackage[margin=1in]{geometry}

\usepackage{hyperref}

\hypersetup{
    pdfstartview={FitH},     % fits the width of the page to the window {FitH},{FitV}
    pdftitle={Estratifica\c{c}\~oes no Espa\c{c}o Moduli dos Fibrados de Higgs},             % title
    pdfauthor={Ronald Alberto Z\'u\~niga-Rojas},     % author
    pdfsubject={Vector bundles on curves and their moduli},           % subject of the document
    pdfkeywords={Moduli of Higgs Bundles, Variations of Hodge Structures, Vector Bundles},           % list of keywords
    pdfnewwindow=true,        % links in new window
    colorlinks=true,          % false: boxed links; true: colored links
    linkcolor=blue,           % color of internal links
    citecolor=red,            % color of links to bibliography
    urlcolor=blue             % color of external links
}

\vfuzz \hfuzz
\newcounter{a}
\ifodd\thea\else\stepcounter{a}\fi

\fancypagestyle{plain}{%
\fancyhf{}

}
\defineshorthand{"-}{\nobreak\hskip0pt\discretionary{-}{-}{-}\nobreak\hskip0pt}

\begin{document}

\theoremstyle{plain}
\newtheorem{Def}{Defini\c{c}\~ao}[section]
\newtheorem{Teo}[Def]{Teorema}
\newtheorem{Cor}[Def]{Corol\'ario}
\newtheorem{Lem}[Def]{Lema}
\newtheorem{Prop}[Def]{Proposi\c{c}\~ao}
\newtheorem{Obs}[Def]{Observa\c{c}\~ao}

\thispagestyle{plain}

\begin{center}
\Large
\textsc{Estratifica\c{c}\~oes no Espa\c{c}o Moduli dos Fibrados de Higgs}
\end{center}

\begin{center}
\bigskip
\normalsize
Julho de 2014\\

\bigskip
\emph{Ronald A. Z\'u\~niga-Rojas\footnote{Este trabalho é financiado por Fundos FEDER através do 
Programa Operacional Factores de Competitividade-COMPETE e por Fundos Nacionais através da FCT-Fundação para a Ciência e a 
Tecnologia no âmbito dos projetos PTDC/MAT-GEO/0675/2012 e PEst-C/MAT/UI0144/2013 e a bolsa de estudo com a refer\^encia 
SFRH/BD/51174/2010.}}\\[6pt]
\small Centro de Matem\'atica da Universidade do Porto \\
\small Faculdade de Ci\^encias da Universidade do Porto \\
\small Rua do Campo Alegre, s/n                    \\
\small 4197-007 Porto, Portugal       \\
\small e-mail: \texttt{ronalbzur@gmail.com}

\end{center}

\noindent
\textbf{Resumo:} 
O trabalho do Hausel prova que a estratifica\c{c}\~ao de Bialynicki-Birula do espa\c{c}o moduli dos fibrados de Higgs de posto dois
coincide com a sua estratifica\c{c}\~ao de Shatz. Estas estratifica\c{c}\~oes n\~ao coincidem para posto geral. Aqui, damos uma abordagem para 
o caso de posto tr\^es da classifica\c{c}\~ao da estratifica\c{c}\~ao de Shatz em termos da estratifica\c{c}\~ao de Bialynicki-Birula. 

\medskip

\noindent
\textbf{Abstract:} 
The work of Hausel proves that the Bialynicki-Birula stratification of the moduli space of rank two Higgs bundles coincides 
with its Shatz stratification. These two stratifications do not coincide in general. Here, we give an approach for the rank three case 
of the classification of the Shatz stratification in terms of the Bialynicki-Birula stratification.

\medskip

\noindent\textbf{palavras-chave:} 
Geometria Alg\'ebrica, Espa\c{c}os Moduli, Teoria de Gauge, Fibrados de Higgs, Estratifica\c{c}\~oes, Fibrados Vetoriais.

\medskip

\noindent\textbf{keywords:} 
Algebraic Geometry, Moduli Spaces, Gauge Theory, Higgs Bundles, Stratifications, Vector Bundles.

\section{Espa\c{c}o Moduli de Fibrados de Higgs}
\label{sec:1}

Seja $\Sigma$\ uma superf\'icie de Riemann compacta de g\'enero $g \geqslant 2$,\ e seja $K = K_{\Sigma} = (T\Sigma)^*$\ o fibrado de linhas 
can\'onico de $\Sigma$.
 \begin{Def}
     Um \textit{fibrado de Higgs} sobre $\Sigma$\ \'e um par $(E, \Phi)$\ onde $E \to \Sigma$\ \'e um fibrado vetorial holomorfo e
     $\Phi: E \to E\otimes K$\ \'e um endomorfismo de $E$\ torcido por $K$,\ chamado \textit{campo de Higgs}. Repare que\ 
     $\Phi \in H^0(\Sigma;\mathrm{End}(E)\otimes K)$.    
 \end{Def} 
Para um fibrado $E \rightarrow \Sigma$\ o \textit{declive} define-se como: $\mu(E):=\frac{\mathrm{deg}(E)}{\mathrm{rk}(E)}=\frac{d}{r}$,\ 
onde $\mathrm{deg}(E)$\ \'e o grau de $E$\ e $\mathrm{rk}(E)$\ \'e seu posto. Para mais informa\c{c}\~ao, veja-se por exemplo {Kobayashi \cite{kob}}.
 \begin{Def}
  Um subfibrado $F \subset E$\ \'e \textit{$\Phi$-invariante} se $\Phi(F)\subset F\otimes K$.\ Um 
  fibrado de Higgs chama-se \textit{semi-est\'avel} se $\mu(F) \leqslant \mu(E)$\ para qualquer subfibrado $\Phi$-invariante 
  n\~ao trivial $0 \neq F\subset E$.\ Chama-se \textit{est\'avel} se a desigualdade \'e estrita para qualquer subfibrado pr\'oprio 
  $\Phi$-invariante n\~ao trivial $0 \neq F\subsetneq E$.\ Finalmente, um fibrado de Higgs chama-se \textit{poli-est\'avel} se \'e 
  a soma directa de subfibrados de Higgs est\'aveis, todos com o mesmo declive.
 \end{Def}
  Fixando o posto $\mathrm{rk}(E)=r$\ e o grau $\mathrm{deg}(E)=d$\ de um fibrado de Higgs $(E,\Phi)$,\ as classes de isomorfismo dos 
  fibrados poli-est\'aveis s\~ao parametrizadas por uma variedade quase-projetiva: o espaço moduli $\mathcal{M}(r,d)$.\ Constru\c{c}\~oes 
  deste espa\c{c}o podem encontrar-se no trabalho de {Hitchin \cite{hit}}, utilizando Teoria de Gauge, ou no trabalho de {Nitsure \cite{nit}}, 
  utilizando m\'etodos de geometria alg\'ebrica.

\section{A\c{c}\~ao de $\mathbb{C}^{*}$\ em $\mathcal{M}(r,d)$}
\label{sec:2}

Existe uma a\c{c}\~ao holomorfa do grupo multiplicativo $\mathbb{C}^{*}$\ em $\mathcal{M}(r,d)$\ definida pela multiplica\c{c}\~ao:
$z\cdot (E,\Phi)\mapsto (E,z\cdot \Phi).$\ Note que Hausel \cite{hau} prova que o limite 
$\displaystyle \lim_{z\to0}z\cdot (E, \Phi) = \lim_{z\to0}(E, z \cdot \Phi)$\ est\'a bem definido e existe para todo o 
$(E,\Phi)\in \mathcal{M}(r,d)$.\ Al\'em disso, este limite \'e fixo pela a\c{c}\~ao de $\mathbb{C}^{*}$.\ Sejam $\{ F_{\lambda} \}$\ as 
componentes irredut\'iveis do lugar de pontos fixos de $\mathbb{C}^{*}$\ em $\mathcal{M}(r,d)$.

Com base no trabalho de Bialynicki-Birula sobre a\c{c}\~oes de grupos alg\'ebricos \cite{b-b}, Hausel \cite{hau} define a estratifica\c{c}\~ao de 
Bialynicki-Birula da seguinte maneira:
  \begin{Def}
    Considere o conjunto $\displaystyle U_{\lambda}^{BB} := \{ (E, \Phi)\in \mathcal{M}(r,d)|\lim_{z\to0}z\cdot (E, \Phi) \in F_{\lambda} \}.$
    Este conjunto \'e o \textit{estrato ascendente} da \textit{estratifica\c{c}\~ao Bialynicki-Birula}:
    $$\mathcal{M} = \bigcup_{\lambda}U_{\lambda}^{BB}.$$
  \end{Def}  
Simpson \cite{sim} prova que os pontos fixos da a\c{c}\~ao $\mathbb{C}^{*} \circlearrowleft \mathcal{M}(r,d)$\ s\~ao as chamadas 
\textit{Varia\c{c}\~oes de Estructura de Hodge}, VHS:
 $$(E,\Phi)\ \textmd{tal que}\ E = \bigoplus_{j=1}^{n}E_j\ \textmd{e}\ \Phi: E_j \to E_{j+1}\otimes K.$$
Dizemos que $(E,\Phi)$\ \'e uma $(\mathrm{rk}(E_1),\ ...,\ \mathrm{rk}(E_n))$-VHS.
   
\section{Estratifica\c{c}\~ao de Shatz}
\label{sec:3}

\begin{Def}\label{HNF}
Uma \textit{filtra\c{c}\~ao Harder-Narasimhan} de $E \to \Sigma$,\ \'e uma filtra\c{c}\~ao da forma:
\begin{equation}
 HNF(E):\ E = E_s \supset E_{s-1} \supset ... \supset E_1 \supset E_0 = 0
\end{equation}
que satisfaz as seguintes propriedades:
\begin{enumerate}[i.]
 \item
 $\mu(E_{j+1}/E_j) < \mu(E_j/E_{j-1})\ \textmd{para}\ 1\leqslant j \leqslant s-1.$
 \item
 $V_j := E_j/E_{j-1}\ \textmd{\'e semi-est\'avel para}\ 1\leqslant j \leqslant s.$
\end{enumerate} 
\end{Def}

\begin{Teo}[{Shatz \cite[Theorem~1]{sha}}]\label{[Sha](Th1)} 
  Todo fibrado vetorial $E \to \Sigma$\ tem uma \'unica filtra\c{c}\~ao de Harder-Narasimhan.
\end{Teo}
Isto foi provado no caso em que $\Sigma$\ \'e uma curva alg\'ebrica projetiva não-singular, por Harder e Narasimhan \cite{hana}. 
A prova de Shatz \cite{sha} \'e v\'alida para variedades projetivas lisas de qualquer dimens\~ao.

\begin{Def}
 Seja $E \to \Sigma$\ um fibrado vetorial de posto $\mathrm{rk}(E)=r$,\ com uma filtra\c{c}\~ao Harder-Narasimhan
 como a mencionada acima em (\ref{HNF}). Definimos o \textit{tipo Harder-Narasimhan}, abreviado HNT, como o vetor
 \begin{equation}
  \mathrm{HNT}(E): \vec{\mu}=(\mu_1,...,\mu_1,\mu_2,...,\mu_2,...\ ...,\mu_s,...,\mu_s) \in \mathbb{Q}^r\label{HNT}
 \end{equation}
 onde $\mu_j = \mu(V_j) = \mu(E_j/E_{j-1})$\ aparece $r_j$ vezes e $r_j=\mathrm{rk}(V_j)$.
\end{Def}

\begin{Def}
    Em consequ\^encia das Proposi\c{c}\~oes 10 e 11 de Shatz \cite{sha}, existe uma estratifica\c{c}\~ao finita de $\mathcal{M}(r,d)$\ 
    pelo tipo Harder-Narasimhan do fibrado vetorial subjacente a um fibrado de Higgs:
    $$\mathcal{M}(r,d) = \bigcup_{t}U'_{t}$$ 
    onde $U'_t \subset \mathcal{M}(r,d)$\ \'e o subespa\c{c}o de fibrados de Higgs $(E,\Phi)$\ cujo fibrado subjacente $E$\ tem 
    $\mathrm{HNT}(E)=t$,\ e a uni\~ao \'e sobre os tipos que existem em $\mathcal{M}(r,d)$.
\end{Def}

 \begin{Prop}[{Hausel \cite[Proposition~4.3.2]{hau}}]\label{[3]Prop4.3.2}
  Se $\mathrm{rk}(E)=2$\ temos que a Estratifica\c{c}\~ao de Shatz coincide com a Estratifica\c{c}\~ao de Bialynicki-Birula.
 \end{Prop} 

\section{Resultado Principal}
\label{sec:4}

Seja $\big[ (E,\Phi) \big] \in \mathcal{M}(3,d)$\ e denote $\displaystyle (E^0,\Phi^0):= \lim_{z\to0}(E,z\cdot \Phi)$. O estrato 
da estratifica\c{c}\~ao Bialynicki-Birula a que $(E,\Phi)$ pertence \'e determinado por $(E^0,\Phi^0)$, e depende do tipo de Harder-Narasimhan 
de $E$, e de certas propriedades de $\Phi$. O nosso Teorema Principal descreve em detalhe esta depend\^encia. 

Para enunciar o Teorema \'e conveniente usar a seguinte notação: para 
um morfismo entre fibrados vetoriais $\phi: E \to F$ vamos escrever $\mathrm{ker}(\phi) \subset E$\ e 
$\mathrm{im}(\phi) \subset F$\ para os subfibrados que se obt\^em saturando os respetivos subfeixes.
\begin{Teo}\label{[ZR](2)}
 \begin{enumerate}[(1.)]
  \item \label{caso1} Suponha que $E\to \Sigma$\ \'e um fibrado holomorfo que tem $\mathrm{HNT}(E) = (\mu_1,\mu_2,\mu_2)$\ 
onde $\mu_j = \mu(V_j)$\ e $V_j = E_j/E_{j-1}$\ s\~ao semi-estáveis. Considere $\phi_{21}: V_1 \to V_2 \otimes K$\ induzido por 
$$E_1 \xrightarrow{\quad \imath \quad} E \xrightarrow{\quad \Phi \quad} E \otimes K \xrightarrow{\jmath \otimes id_K} \big( E/E_1\big) \otimes K.$$
Defina $\mathcal{I}:=\phi_{21}(E_1)\otimes K^{-1} \subset V_2$\ onde $\mathrm{rk}(\mathcal{I})=1$,\ e defina tamb\'em 
$F:=V_1\oplus \mathcal{I} \subset V_1 \oplus V_2 = E$\ onde $\mathrm{rk}(F)=2$.\ Ent\~ao temos duas possibilidades:
    \begin{enumerate}
      \item[(1.1.)]
      Suponha que $\mu(F)<\mu(E)$.\ Ent\~ao, $(E^0,\Phi^0)$ \'e uma $(1,2)$-VHS da forma:
      $$(E^0,\Phi^0) = \Big(V_1\oplus V_2, \left( \begin{array}{c c}
                            0     & 0\\ 
                        \phi_{21} & 0
                       \end{array}
      \right) \Big).$$
      \item[(1.2.)]
      Por outro lado, se $\mu(F) \geqslant \mu(E)$,\ ent\~ao, $(E^0,\Phi^0)$ \'e uma $(1,1,1)$-VHS da forma:
      $$(E^0,\Phi^0) = \Big(L_1\oplus L_2\oplus L_3, \left( \begin{array}{c c c}
                            0        & 0 & 0\\ 
                        \varphi_{21} & 0 & 0\\
                            0        & \varphi_{32} & 0
                       \end{array}
      \right) \Big)$$
      onde $L_1, L_2, \textrm{and}\ L_3$\ s\~ao fibrados de linhas.
    \end{enumerate}

  \item \label{caso2} Analogamente, suponha que $E\to \Sigma$\ \'e um fibrado holomorfo tal que $\mathrm{HNT}(E) = (\mu_1,\mu_1,\mu_2)$\ 
  onde $\mu_j = \mu(V_j)$\ e $V_j = E_j/E_{j-1}$\ s\~ao semi-estáveis. Considere $\phi_{21}: V_1 \to V_2 \otimes K$\ induzido por
$$E_1 \xrightarrow{\quad \imath \quad} E \xrightarrow{\quad \Phi \quad} E \otimes K \xrightarrow{\jmath \otimes id_K} \big( E/E_1\big) \otimes K.$$
Defina $N:=\mathrm{ker}(\phi_{21})\subset V_1$\ onde $\mathrm{rk}(N)=1$.\ Ent\~ao, temos duas possibilidades:
    \begin{enumerate}
      \item[(2.1.)]
      Suponha que $\mu(N)<\mu(E)$.\ Ent\~ao, $(E^0,\Phi^0)$\ \'e uma $(2,1)$-VHS da forma:
      $$(E^0,\Phi^0) = \Big(V_1\oplus V_2, \left( \begin{array}{c c}
                            0     & 0\\ 
                        \phi_{21} & 0
                       \end{array}
      \right) \Big).$$      
      \item[(2.2.)]
      Por outro lado, se $\mu(N) \geqslant \mu(E)$,\ ent\~ao: $(E^0,\Phi^0)$\ \'e uma $(1,1,1)$-VHS da forma:
      $$(E^0,\Phi^0) = \Big(L_1\oplus L_2\oplus L_3, \left( \begin{array}{c c c}
                            0        & 0 & 0\\ 
                        \varphi_{21} & 0 & 0\\
                            0        & \varphi_{32} & 0
                       \end{array}
      \right) \Big)$$
      onde $L_1, L_2, \textrm{and}\ L_3$\ s\~ao fibrados de linhas.
    \end{enumerate}

  \item \label{caso3} Finalmente, suponha que $E\to \Sigma$\ \'e um fibrado holomorfo tal que $\mathrm{HNT}(E) = (\mu_1,\mu_2,\mu_3)$\ 
  onde $\mu_j = \mu(V_j)$\ e $V_j = E_j/E_{j-1}$\ s\~ao semi-estáveis.
  \begin{enumerate}
   \item[(3.1.)]
   Suponha que $\mu(E_2/E_1)<\mu(E)$.\ Ent\~ao podemos definir $F$\ como no caso (\ref{caso1}.), e ent\~ao, 
temos duas possibilidades:
      \begin{enumerate}
	\item[(3.1.1.)] Suponha que $\mu(F)<\mu(E)$.\ Ent\~ao: $(E^0,\Phi^0)$ \'e uma $(1,2)$-VHS.
	\item[(3.1.2.)] Se $\mu(F) \geqslant \mu(E)$,\ ent\~ao: $(E^0,\Phi^0)$ \'e uma $(1,1,1)$-VHS.
      \end{enumerate}

   \item[(3.2.)]
Por outro lado, se $\mu(E_2/E_1)>\mu(E)$,\ ent\~ao podemos definir $N$\ como no caso (\ref{caso2}.), e ent\~ao, 
temos duas possibilidades:
    \begin{enumerate}
      \item[(3.2.1.)] Se $\mu(N)<\mu(E)$,\ ent\~ao: $(E^0,\Phi^0)$\ \'e uma $(2,1)$-VHS.
      \item[(3.2.2.)] Se $\mu(N) \geqslant \mu(E)$,\ ent\~ao: $(E^0,\Phi^0)$\ \'e uma $(1,1,1)$-VHS.
    \end{enumerate}   
  \end{enumerate}
 \end{enumerate}
\end{Teo}

\newpage 

\renewcommand*{\refname}{Refer\^encias}

\begin{flushright}
  \textbf{Updated filiation and address:}\\
  \emph{Ronald A. Z\'u\~niga-Rojas}\\
  \small Centro de Investigaciones Matem\'aticas\\
  \small y Metamatem\'aticas CIMM\\
  \small Universidad de Costa Rica UCR\\
  \small San Jos\'e 11501, Costa Rica\\
  \small e-mail: \texttt{ronald.zunigarojas@ucr.ac.cr}\\
  \small $2017-10-26$.
\end{flushright}

\end{document}